\def\year{2019}\relax
\documentclass[letterpaper]{article} 
\usepackage{aaai19}  
\usepackage{times}  
\usepackage{helvet}  
\usepackage{courier}  
\usepackage{url}  
\usepackage{graphicx}  
\usepackage{bm,amsfonts,amssymb,amsmath,dsfont}
\usepackage{multirow,booktabs}

\newtheorem{prop}{Proposition}
\usepackage{algorithm}
\usepackage[noend]{algpseudocode}
\usepackage{algpseudocode}
\algnewcommand{\Inputs}[1]{%
  \State \textbf{Inputs:}
  \Statex \hspace*{\algorithmicindent}\parbox[t]{.8\linewidth}{\raggedright #1}
}
\algnewcommand{\Initialize}[1]{%
  \State \textbf{Initialize:}
  \Statex \hspace*{\algorithmicindent}\parbox[t]{.8\linewidth}{\raggedright #1}
}

\algnewcommand{\Outputs}[1]{%
  \State \textbf{Outputs:}
  \Statex \hspace*{\algorithmicindent}\parbox[t]{.8\linewidth}{\raggedright #1}
}

\begin{document}
%
\title{The Contextual Appointment Scheduling Problem}
\author{Nima Salehi Sadghiani \\
Department of Industrial and Operations Engineering\\
University of Michigan\\
Ann Arbor, MI 48105\\
\texttt{nsalehi@umich.edu}
\And
Saeid Motiian \\
Adobe Inc.\\
San Jose, California, USA\\
\texttt{motiian@adobe.com}}

\maketitle
\begin{abstract}
This study is concerned with the determination of optimal appointment times for a sequence of jobs with uncertain duration. We investigate the data-driven Appointment Scheduling Problem (ASP) when one has $n$ observations of $p$ features (covariates) related to the jobs as well as historical data. We formulate ASP as an Integrated Estimation and Optimization problem using a task-based loss function. We justify the use of contexts by showing that not including them yields to inconsistent decisions, which translates to sub-optimal appointments. We validate our approach through two numerical experiments.
\end{abstract}

\section{Introduction}  \label{section1}
Appointment Scheduling Problem (ASP) arises in many service industries to improve the service quality, operational efficiency, utilization of resources, match workload to available capacity, and smooth the flow of customers (e.g. health care, legal services, accounting services, and loading/unloading ships/planes). In health care industry, ASP emerges in numerous settings, such as scheduling outpatient appointments in primary care \cite{cayirli2006designing}, surgery scheduling \cite{denton2003sequential}, and call-center/nurse staffing \cite{gurvich2010staffing,koccauga2015staffing}.

Assuming punctuality, a common problem faced by decision makers is how to determine the starting time of job when its duration is uncertain. The decision makers' goal is to minimize the expected costs associated with jobs' waiting time, server idle time, and overtime. For example, when appointments are set, jobs are not available prior to the scheduled start time, even if the server becomes idle. Thus, at the cost of additional waiting (for jobs), choosing an early start time will lead to better server utilization, whereas a late start time will reduce the waiting time at the cost of additional server idling.

A common assumption in the literature of the ASP is that we can accurately estimate the probability distribution of the job duration. However, in reality, this may not be the case. The scarcity of the historical data, diversity in types of jobs can easily violate the underlying assumptions behind distribution fitting (like independent and identically distributed (i.i.d.) assumption, the choice of distribution, or the correlation among the job duration distributions). Although the ASP is a well-explored problem, most studies in the literature made these assumptions without considering the consequences.

In this paper, we investigate the ASP under two settings: classic models in which we only have access the previous job duration, and the feature-based models in which we have access to some exogenous information about the job duration as well as the actual value. Rather than a two-step process of first estimating the duration distribution then optimizing for the optimal scheduling, we propose solving the ASP via a single step by training a Deep Neural Network (DNN) in an end-to-end task-based framework. The difficulty of the integrated learning and optimizing task in ASP arises from the fact that the decision about each job duration depends on the actual realization of the previous jobs.

\section{Related Work}  \label{section2}
\subsection{ASP literature}
ASP studies with uncertain service times often assume full distributional information about the job durations (i.e. full joint distribution of the jobs is known, or jobs are independent and marginal distributions are known). Excellent surveys of ASP models can be found in\cite{gupta2008appointment,cardoen2010operating}.

Sample Average Approximation (SAA) is a common method used to solve the stochastic ASP. SAA approach can be computationally intensive, and require knowledge of the joint distribution of the job service durations \cite{denton2003sequential}.

Further, when the joint distribution of the uncertain parameters is not known to us, other types of modeling is proposed in the literature. \citeauthor{mittal2014robust} \citeyear{mittal2014robust} proposed Robust Optimization (RO) for the ASP where job  durations lie within interval uncertainty sets. \citeauthor{mak2014appointment} \citeyear{mak2014appointment} proposed Distibutionally Robust Optimization (DRO) for ASP assuming partial distributional information of the uncertain parameters such as support and moments are known. The appointments are made to minimize the worst-case expected objective value among all possible distributions with the specified support and moments information. Finally, \citeauthor{zhang2017distributionally} \citeyear{zhang2017distributionally} proposed using DRO with moment uncertainty set and derived an approximate semidefinite programming model.

\subsection{Prescriptive Analytics}
Solving many real world problems requires two main steps of estimation (prediction) and optimization (prescription). Due to the complex nature of each step, researchers typically solve these two steps separately. The major issue with this approach is that the optimization problem is blind to the prediction process and its sensitivity to over/under estimation of the parameter. In contrast to the Separated Estimation and Optimization (SEO) approach, another framework is proposed by \citeauthor{bertsimas2014predictive} \citeyear{bertsimas2014predictive}. They referred to this method as "prescriptive analytics". They combined the ideas from Machine Learning (ML) and Operations Research and Management Science (OR/MS) in developing a framework for using auxiliary data to prescribe optimal decisions.

A general framework for the Integrated Estimation and Optimization (IEO) approach was proposed later by \citeauthor{elmachtoub2017smart} \citeyear{elmachtoub2017smart}. They called their framework "Smart Predict, then Optimize” (SPO). Training a model with respect to the SPO loss is computationally challenging, and therefore they also proposed a surrogate loss function, called the SPO+ loss, which upper bounds the SPO loss, has desirable convexity properties, and is statistically consistent under mild conditions. Minimizing an application-specific loss function is refereed to as "training with direct loss" \cite{song2016training} or task-based model learning \cite{donti2017task} in the ML literature. In this paper, we proposed using DNN to solve the feature-based ASP.

\section{ASP Modeling}  \label{section3}
We consider an ASP with a single server at which jobs arrive
punctually at scheduled appointment times. The set of jobs follows a fixed order of arrivals. We formulate the stochastic ASP with random service duration. Throughout this paper, we use boldface notation to denote vectors.

\subsection{ASP Classic Models}  \label{section31}
We consider $n$ jobs arriving at a single server given a predetermined sequence $1,2,3,\dots,n$. Job $i$ has a random service duration ${p_{i}}$. The planner determines time allowances $S_{i} \geq 0$ for each job $i$. If job $i$ cannot be started at its planned start time due to a delay of completion of the previous job, a waiting time cost will be incurred per unit time of delay. Let $W_{i}$ denote the waiting time of the job $i$. On the other hand, if job $i-1$ completed sooner than its planned start time, an idling time cost will be incurred per time. Let $I_{i}$ denote the idling time of the server between job $i-1$ and job $i$. The goal is to optimize the decisions of scheduling for arrival times of each job,
or equivalently, assigning time intervals between each two jobs.

The random job durations $\bm{p}$ follow a joint probability distribution $F(\bm{p})$. For any job $i$, the support of ${p_{i}}$ is independent from other jobs and is denoted by $\mathds{D}_i$; therefore, $\mathds{D} = \mathds{D}_1 \times \dots \times \mathds{D}_n$. Given time allowances $\bm{S} \in \bm{\Re^{+}}$ and a realization of the random service times, denoted by $\bm{p}$, the waiting and idling times should satisfy the following constraints:

\begin{align}
    & W_{i} - I_{i} = W_{i-1}+p_{i-1}-S_{i-1} & \forall2 \leq i \leq n \label{baseConst}
 \end{align}

Moreover, $W_{i}, I_{i} \in \Re^{+}$. It is also worth noting that because of the punctually arrival assumption $W_{1} = I_{1} =0$. Assuming linear costs for waiting and idling, the total waiting and idling cost is denoted by
\begin{align}
    & f(\bm{S},\bm{p}) = c^{\mathcal{W}}\sum_{i=1}^{n}W_{i}+c^{\mathcal{I}}\sum_{i=1}^{n}I_{i} \label{baseObj}
\end{align}

where $c^{\mathcal{W}}$ and $c^{\mathcal{I}}$ (non-negative) represent the costs of waiting and idling, respectively.

\subsubsection{ASP Model with Distributional Information.}  \label{section311}
Assuming the joint probability distribution of job durations $F(\bm{p})$ is known, we can formulate the stochastic ASP model as follow
\begin{align}
    & \min_{\bm{S} \in \Re^{+}} \mathds{E}_{\bm{p} \sim F(.)}\mathcal{F}(\bm{S},\bm{{p}}) \label{expObj}
\end{align}

where $\mathds{E}$ denotes the expected value and $\mathcal{F}(.) := \min_{\bm{W}, \bm{I} \in \Re^{+} \cup \{eq. (\ref{baseConst})\}} f(.)$. The objective function of the stochastic ASP is to minimize the expected total waiting and idling costs under the assumption that random job durations come from the known joint probability distribution $F(.)$ with support $\mathds{D}$.

\subsubsection{ASP Model with Limited Distributional Information.} \label{droASP}
A key assumption of the previous model is that the exact joint
distribution $F(.)$ is known to the planner at the beginning. However, this assumption may not be held under different practical settings. To address the distributional ambiguity issue, in this section, we assume only some limited distributional information about $F(.)$ is available to the planner as well as the support $\mathds{D}$. We consider a moment-based ambiguity set for the distributions of the uncertain parameters.

For moment-based ambiguity set, the $q$th moment of ${p}_{i}$ is denoted by $M_{iq}$. We also use $\mu_{i}$ and $\sigma_i$ to denote the mean and standard deviation of ${p}_{i}$, respectively; that is, $M_{i1} = \mu_i$ and $M_{i2} = \mu^{2}_{i}+\sigma^{2}_{i}$. For any $q \in Q$ where $Q$ is a finite set of positive integers, the distribution $F(.)$ belongs to the ambiguity set $\mathcal{M}(\mathds{D},Q)$
\begin{align}
  & :=
  \left \{ \begin{aligned}
	& \int_{\mathds{D}}dF(\bm{p})=1 &\\
    & \int_{\mathds{D}}p_{i}^{q}dF(\bm{p})=M_{iq} & \forall q \in Q,  \qquad  1 \leq i \leq n & \\
    & dF(\bm{p}) \geq 0 & \forall \bm{p} \in \mathds{D}
    \end{aligned}
  \right \} \label{momAmb}
\end{align} 

Any set of $\bm{p}$ belong to $\mathcal{M}(.)$ is a distribution, match the moments, and it is non-negative. Under this assumption, we can formulate the stochastic ASP as a min–max problem robust to the distribution ambiguity
\begin{align}
    & \min_{\bm{S} \in \Re^{+}} \max_{F \in \mathcal{M}} \mathds{E}_{\bm{p} \sim F(.)}\mathcal{F}(\bm{S},\bm{{p}}) \label{minmaxObj}
\end{align}

We shall refer to this framework as Distributionally Robust Optimization (DRO).

\subsubsection{Data-driven ASP without Distributional Assumptions.}
In practice, the planner does not know the true joint distribution of the job durations. If one has access to historical job durations for $T$ periods: $(\bm{p}^t) = (p^t_{1},\dots,p^t_{n})$, $1 \leq t \leq T$, but no other information (no context), then the sensible approach is to substitute the true expectation with a sample average expectation and solve the resulting problem.
\begin{align}
    & \min_{\bm{S} \in \Re^{+}} \frac{1}{T} \sum_{t \in T}\mathcal{F}(\bm{S},\bm{p}^t) \label{saaObj}
\end{align}

This approach is called the Sample Average Approximation (SAA) approach in stochastic optimization \cite{kleywegt2002sample}. The consistency of the SAA estimator for the expectation of the objective function of the stochastic models has been discussed thoroughly in \cite{shapiro2009lectures}. 

\subsection{Drawbacks of the Classic Models}
In the previous section, we reviewed classic approaches to model ASP. In a realistic situation, the classic models are too simplistic to represent many real scenarios because one can collect data on exogenous information about the jobs' characteristics as well as the job durations. In other words, the planner has access to a richer data set for decision making. Further, most of the models fail to consider context in estimating the job durations. The classic models obtain an estimate of the probability distribution and solve the optimization problem separately. The approximation process of the probability distributions involves errors, especially when historical data is scarce or when durations are not i.i.d.

\subsection{The Feature-based ASP Models}
In practice, the job durations depend on many observable features (equivalently, covariates or context). For example, in patient appointment scheduling, the characteristics of the patients can be used as features. If the jobs are flight arrivals, then the attributes of the flights can be used (i.e. day, month, season, weather, origin and destination airports). 

Assuming the planner schedules $n$ jobs with $m$ features at each time period and s/he has access to historical data for $T$ periods, we have $\{(\bm{x}^t,\bm{p}^t)\}_{t=1}^T$ where $\bm{x}^t = (x_1^t, \dots, x_n^t)$ and $\bm{p}^t = (p_1^t, \dots, p_n^t)$. For notation brevity, we denote all the $m$ features for the job $i$ at time $t$ with $x_i^t, 1 \leq i \leq n$, i.e. $x_i^t = (x_{i1}^t, \dots, x_{im}^t)$.

In the feature-based ASP model, the features are available to the planner prior to the decision making process. In other words, the decision maker optimizes the conditional expected cost function:
\begin{align}
    & \min_{\bm{S}(.) \in \mathcal{S},\{\bm{S}:  \mathcal{X} \rightarrow  \Re^{+}\} } \mathds{E}_{t}f(\bm{S}(\bm{x}^t),\bm{p}^t|\bm{x}^t) \label{condObj}
\end{align}

where the scheduling decisions $\bm{S}(.)$ are now functions mapping the feature space $\mathcal{X} \subset \Re^m$ to the positive real numbers ($\Re^{+}$) and the expected waiting and idling costs that we minimize is now conditional on the feature vector $\bm{x}^t \in \mathcal{X}$.

The planner needs to be able to solve the described problem in equation (\ref{condObj}) efficiently given the observation $\bm{x}^t$. The main issue in this model is that the job durations might be correlated and allowance decisions have to be made sequentially. For example, when the first job is being processed on the server, we do not know whether the server will be idle until the next job arrives or the next job will be delayed. The uncertainty about the starting time of each job propagates through the sample of $n$ jobs at each time unit $t$. The planner needs to know the actual realization of the previous job and its starting time to optimally plan for the next one.

\subsubsection{2-Job ASP Model.}
The special case of this problem is when we have $n = 2$. Because of punctual arrivals the first job arrives at time zero. For this case, it is only necessary to predict the duration of the first job, and the second job starts its process right when the first job is done. \citeauthor{weiss1990models} \shortcite{weiss1990models} showed this special case corresponds to the \textit{Newsvendor} problem. Therefore, under the assumption that no data features is available but the underlying distribution of the job durations is known, a closed form expression for the optimal allowance of the first job can be obtained. The closed form solution for the classic \textit{Newsvendor} problem described in \citeauthor{gallego1993distribution} \shortcite{gallego1993distribution} is as follow
\begin{align}
    \hat{p}^{*}_1 = S^{*}_1 = F^{-1}\left(\frac{c^{\mathcal{I}}}{c^{\mathcal{I}}+c^{\mathcal{W}}}\right) \label{nvSol}
\end{align}

where $\hat{p}^{*}_1$ is the optimal estimation for the duration of the first job and also its time allowance. The $F(.)^{-1}$ refers to the inverse distribution of $F(.)$. Recently, some scholars used ML techniques to solve the data-driven \textit{Newsvendor} problem. Here, we discuss the ones that are closely related to our model.

\citeauthor{ban2017big} \shortcite{ban2017big} has shown analytically the superiority of using Linear Model (LM) to describe the optimal solution than other data-driven approaches such as SAA. Basically, in their models they defined the mapping function $\bm{S}(\bm{x}^t)$ as:
\begin{align}
    & \mathcal{S}_{LM} = \left\{\bm{S}:  \mathcal{X} \rightarrow \Re, \bm{S}(\bm{x}^t) = \sum_{k=0}^{m}\beta_{k}x^{t}_{k}\right\}  \label{lmSol}
\end{align}

where $\beta_k$s are the coefficients of the linear model and the $x^t_0=1$ accommodates for a feature-independent term. Further, they have shown that a complex-valued $\bm{S}(\bm{x}^t)$ mapping functions (but infinitely differentiable at $0$) can be approximated by polynomial terms by Taylor series. In their setting, the objective function becomes:
\begin{align}
    & \min_{\bm{S}(.) \in \mathcal{S}_{LM},\{\bm{S}: \mathcal{X} \rightarrow \Re\}} \frac{1}{T}\sum_{t=1}^{T}f\left(S_1(\bm{x}^t),p_1^t\right)  = \nonumber\\
    & \qquad \frac{c^{\mathcal{W}}}{T} \left \lVert \left(\bm{p}_1-S_1(\bm{x})\right)_{+}\right \rVert_{1}+\frac{c^{\mathcal{I}}}{T} \left \lVert \left (S_1(\bm{x})-\bm{p}_1\right)_{+}\right \rVert_{1} \label{lmObj}
\end{align}

In equation (\ref{lmObj}), the objective function has the form of $\mathcal{L}_1$ asymmetric loss function where $c^{\mathcal{W}}$ and $c^{\mathcal{I}}$ represent the unit cost of overestimation and underestimation, respectively. This special structure is due to that fact that $W_1^t \cdot I_1^t=0, 1 \leq t \leq T$. Further, they discuss the conditions for the optimality of the decisions from the feature-based model and the asymptotics of the optimal solutions. Finally, they have shown the decisions from SAA (with finite sample sizes) has bias of $\mathcal{O}(1)$; hence the SAA decisions may not be asymptotically optimal.

\citeauthor{OroojlooyjadidS16} \shortcite{OroojlooyjadidS16} have shown the benefits of using non-linear models to estimate the mapping function $\bm{S}(\bm{x}^t)$ in comparison with linear models. In their model $\bm{S}(\bm{x}^t)$ is a complex and non-linear function (dense layer architecture).
\begin{align}
    & \mathcal{S}_{DNN} = \left\{\bm{S}:  \mathcal{X} \rightarrow \Re, \bm{S}(\bm{x}^t,\bm{\omega},\bm{b}) \right\} \label{dnnSol}
\end{align}

where $\bm{\omega}$, $\bm{b}$ are the weights and biases for a fully-connected network architecture. They used sigmoid functions to create non-linearity in the model. They also suggested using a surrogate $\mathcal{L}_2$ loss function for faster convergence.
{\small
\begin{align}
    &\min_{\bm{S}(.) \in \mathcal{S}_{DNN},\{\bm{S}: \mathcal{X} \rightarrow \Re\}} \frac{1}{2T}\sum_{t=1}^{T}f^{2}\left(S_1(\bm{x}^t,\bm{\omega},\bm{b}) ,p_1^t\right) = \nonumber\\
& \frac{c^{\mathcal{W}}}{2T}\left \lVert \left(\bm{p}_1-S_1(\bm{x},\bm{\omega},\bm{b})\right)_{+}\right \rVert^{2}_{2}+ \frac{c^{\mathcal{I}}}{2T} \left \lVert \left (S_1(\bm{x},\bm{\omega},\bm{b})-\bm{p}_1\right)_{+}\right \rVert^{2}_{2} \label{dnnObj}
\end{align}
}%

\citeauthor{zhang2017assessing} \shortcite{zhang2017assessing} discuss the benefits of using an extra layer of Rectified Linear Units (ReLUs) to decompose the original $\mathcal{L}_1$ loss function into a constant ($p^t_1$) and a differentiable function ($S_1(\bm{x}^t,\bm{\omega},\bm{b})$). They used $c^{\mathcal{W}}$ and $c^{\mathcal{I}}$ as the weights of the new layer to show the original $\mathcal{L}_1$ loss function can be reconstructed from these two parts.

\subsubsection{$n$-Job ASP Model.}
Assuming $n$ jobs are needed to be scheduled, for scheduling each job, the planner needs to take into account the actual realizations of the previous jobs. One can clearly see how the previous methods in the literature fail to consider this recursive relationship among the predictors (decision variables of the optimization model). To address this shortcoming, we propose a different loss function to learn the complex and recursive behavior of the job durations from the context while it accounts for the sequential relationships among jobs that are required to be scheduled in each time unit.

To estimate the mapping function $\bm{S}(\bm{x}^t)$, we propose adding a sigmoid layer on the outputs of the dense layer to satisfy the non-negativity of the decision variables. We will refer to this mapping as $DNN^+$,
\begin{align}
    & \mathcal{S}_{DNN^+} = \left\{\bm{S}:  \mathcal{X} \rightarrow \Re^{+}, \bm{S}(\bm{x}^t,\bm{\omega},\bm{b}) \right\} \label{dnn+Obj}
\end{align}

Further, we suggest using Quantile Regression (QR) Loss Function for the neural network. \citeauthor{white1992nonparametric} \shortcite{white1992nonparametric} presents theoretical support for the use of QR for an ANN to estimate the potentially non-linear quantile models. \citeauthor{taylor2000quantile} \shortcite{taylor2000quantile} has shown the application of QR-NN to estimate the conditional density of the uncertain parameters in a multi-period forecasting framework. Minimizing the Mean Absolute Deviation (MAD) leads to an estimate of the conditional median of the predict and data. By applying asymmetric weights to positive/negative errors by using a tilted form of the absolute value function, one can instead compute conditional quantiles of the predictive distribution \cite{koenker1978regression}. The tilted absolute value function (also known as the pinball loss function) for each individual prediction is given by
\begin{align}
   & \mathcal{L}(p_i-\hat{p}_i|q) = \begin{cases}
    q(p_i-\hat{p}_i), & \text{if } p_i-\hat{p}_i \geq 0,\\
    (1-q)(\hat{p}_i-p_i), & \text{otherwise}.
  \end{cases} \label{pinball}
\end{align}

where $q$ is the quantile of interest. We can show the loss function for the ASP can be represented as follows:\\

\begin{prop}
The proposed surrogate loss function for the ASP can be represented by
{\small
\begin{align}
& \min_{\bm{S}(.) \in \mathcal{S}_{DNN^+},\{\bm{S}: \mathcal{X} \rightarrow \Re^{+}\}} \frac{1}{2T(c^{\mathcal{W}} + c^{\mathcal{I}})}\sum_{t=1}^{T}f\left(\hat{p}^t ,p^t\right) = \nonumber\\
& \frac{1}{2T(c^{\mathcal{W}} + c^{\mathcal{I}})}\left(c^{\mathcal{W}}\left \lVert \sum_{i=2}^{n} \left(\bm{W}_{i-1} + \bm{p}_{i-1}-\bm{\hat{p}}_{i-1}\right)_{+}\right \rVert_{1} \right)+ \nonumber\\
&  \qquad \qquad \qquad \left( c^{\mathcal{I}} \left \lVert \sum_{i=2}^{n} \left (-\bm{W}_{i-1}+\bm{\hat{p}}_i-\bm{p}_i\right)_{+}\right \rVert_{1}\right) \label{ASPloss}
\end{align}}

where $q = \frac{c^{\mathcal{W}}}{c^{\mathcal{W}} + c^{\mathcal{I}}}$, $\hat{p}_{i} = \bm{S}(\bm{x}_{i},\bm{\omega},\bm{b})$, and $W_{i}$ is calculated from a modified version of Lindley’s recursion \shortcite{lindley1952theory}. Note that we use $(.)_{+}$ indicating $\max(. , 0)$. 
\end{prop}

\section{Optimization Methods for ASP Models}  \label{section4}
\subsection{Approaches for Classic ASP}
For the classic ASP models, we propose using SAA for the stochastic optimization models. This approach is applicable when the joint distribution of the job durations is known or when we have enough historical data. In the former case, we sample from the distribution and generate scenarios. In the latter, we treat each historical realization as a scenario. If limited distributional information is available, we propose using a cutting-plane approach to solve the DRO models.

\subsubsection{SAA for Two-stage Stochastic Optimization.}
Having known the distribution of the job durations, we can sample and use the SAA scheme to solve the two-stage stochastic optimization model. Similar to  \citeauthor{denton2003sequential} \shortcite{denton2003sequential}, we can rewrite the ASP problem as follow into a two-stage optimization model:
{\small
\begin{subequations} \label{saa}
\begin{align}
    & \min_{\bm{S},\bm{W},\bm{I} \in\Re^{+}} \frac{1}{N} \sum_{j=1}^{N}\left(c^{\mathcal{W}}\sum_{i=1}^{n}W_{i}(\xi_j)+c^{\mathcal{I}}\sum_{i=1}^{N}I_{i}(\xi_j)\right) \label{saaobj}\\
    \nonumber\\
    & W_{i}(\xi_j) \geq W_{i-1}(\xi_j)+p_{i-1}(\xi_j)-S_{i-1}, \forall \xi_j \in \Omega, 2 \leq i \leq n\\
    &  I_{i}(\xi_j) \geq -W_{i-1}(\xi_j)-p_{i-1}(\xi_j)+S_{i-1}, \forall \xi_j \in \Omega,  2 \leq i \leq n
\end{align}
\end{subequations}
}%
where $\xi_j$ is a scenario in scenario set $\Omega$. Each scenario is a set of realizations for job durations $1 \leq i \leq n$.

\subsubsection{A Cutting-plane algorithm for DRO.}
When we have limited distributional information, we propose modeling ASP as a DRO model. \citeauthor{mak2014appointment} \shortcite{mak2014appointment}  formulated the equivalent models for the min-max problem with moment-based ambiguity set under some specific assumptions. In here, we develop a generic cutting-plane algorithm to solve the min-max model for the moment-based ambiguity set. The generic form of the ASP DRO model is as follow
\begin{align}
    & \min_{\bm{S} \in \Re^{+}} \max_{F \in \mathcal{M}} \mathds{E}_{\bm{p} \sim F(.)}\mathcal{F}(\bm{S},\bm{{p}}) \nonumber\\
        & \qquad \equiv \min_{\bm{S} \in \Re^{+}} \max_{F \in \mathcal{M}} \int_{\mathds{D}}\mathcal{F}(\bm{S},\bm{{p}})dF(\bm{p}) \label{reformObj}
\end{align}

Although this model is computationally intractable, an approximation of this problem with discrete job durations can be solved with the following cutting-plane algorithm
\begin{algorithm}
\caption{A cutting-plane algorithm for DRO ASP}\label{benders}
\begin{algorithmic}[1]
    \State $\text{Set tolerance level } \epsilon \gets (0,1)$
    \State $\text{Initialize iteration counter } k \gets 0$
    \State $\text{Initialize cuts list } cuts \gets \emptyset$
        \State $\text{Set } g(\bm{S}^{0})  \gets \infty$
    \State $\text{Set stop condition } newCuts \gets \text{True}$
    \While{$newCuts$}
        \State $k \gets k+1$
        \State Solve \textit{(Master)}:
        \begin{align*}
        & \min_{\bm{S},\bm{\alpha}, \delta} \sum_{i=1}^{n}\sum_{q \in Q}\alpha_{iq}M_{iq}+\delta \\
            & \text{subject to:} \notag\\
            & \qquad \delta \geq g(\bm{S}^k), \quad k=0,1,\ldots,K
        \end{align*}
        \State Let $\bm{\bar{S}}^k, \bm{\bar{\alpha}}^k$, and $\bar{\delta}^k$ be the optimal solutions for the \textit{(Master)}, then solve \textit{(Subproblem)}:
        \begin{align*}
        &  g(\bm{\bar{S}}^k) \equiv \max_{\bm{\hat{p}} \in \{0,1\}^{n*L},\bm{\lambda}\in \Re^{+},\bm{Y} \in \digamma^D(\bm{Y})}\\
        & \qquad \sum_{i=1}^{n}\sum_{l \in L_i} l \lambda_{il} - \sum_{i=1}^{n}\bar{S}^k_{i-1}Y_i - \sum_{i=1}^{n}\sum_{q \in Q}\sum_{l \in L_i}\bar{\alpha}^{k}_{iq}l^q\hat{p}_{il}\\    
        & \text{subject to:} \notag\\
        & \sum_{l \in L_i} \hat{p}_{il}=1, \forall i\\
   &\lambda_{il} \geq -c^{\mathcal{I}} \hat{p}_{(i-1)l}\\
   &\lambda_{il} \geq Y_{i} - M(1-\hat{p}_{(i-1)l})\\
   &\lambda_{il} \leq M\hat{p}_{(i-1)l}\\
   &\lambda_{il} \leq Y_{i} + c^{\mathcal{I}}(1-\hat{p}_{(i-1)l}) 
        \end{align*}
         \If{$\bar{\delta}^k < (1-\epsilon)g(\bm{\bar{S}}^k)$}
            \State $cuts \gets cuts \cup \{ \delta^{k} \geq g(\bm{S}^k)\}$
        \Else
            \State $newCuts \gets \text{False}$
        \EndIf
    \EndWhile
\end{algorithmic}
\end{algorithm}

\subsection{Approaches for Feature-based ASP}
In the feature-based ASP models, we incorporate the features to predict the job durations. This can be done in two settings: (1) SEO: Separate Estimation and Optimization, (2) IEO: Integrate Estimation and Optimization.

\subsubsection{Separate Estimation and Optimization (SEO).}
In this approach, we incorporate the data feature information in the ASP by first mapping the job durations on the feature space assuming a normally distributed error term; then, we treat the estimates as the deterministic values in the ASP model. \citeauthor{ban2017big} \shortcite{ban2017big} refer to this approach as SEO. In this paper, we use neural nets to estimate the job durations from the feature space. Then, we solve the optimization problem given the predicts from the network.

\subsubsection{Integrate Estimation and Optimization (IEO) with QR-NN.}
In this approach, we try to find a parameterization of the features that optimizes the objective function of the ASP. In other words, we try to find a mapping function from the feature space to the optimal appointment decisions. To do so, we minimize the surrogate loss function in equation (\ref{ASPloss}) (direct loss for the neural network). The estimates of this loss function are the optimal decisions of the stochastic ASP. In general, computing the gradients of an objective that depends upon the argmin/argmax operation is challenging. Recently, several authors developed different techniques for argmin/argmax differentiating \cite{gould2016differentiating,donti2017task}. For the ASP, due to the special structure of the optimization problem, we do not need to worry about this issue. However, still the gradients of the loss (\ref{ASPloss}) are convolutional and hard to derive. Specifically, for the loss function $\mathcal{L}(\bm{\theta})$:
\begin{align}
    & \frac{\partial \mathcal{L}(\bm{\theta})}{\partial \bm{\theta}} = \sum_{i=2}^{n}\frac{\partial \mathcal{L}(\bm{\theta})}{\partial \bm{\hat{p}}_i} \frac{\partial  \bm{\hat{p}}_i(\bm{\theta})}{\partial \bm{\theta}} \qquad \qquad  \label{derv}
\end{align}

The first term in equation (\ref{derv}) is the partial derivative of the loss function to the optimal job duration $i$ and the second term is the partial derivative of each estimate with respect to the
network parameters $\bm{\theta}$ (weights and biases).  Calculating both terms in equation (\ref{derv}) is a complicated task since it requires the estimate and actual job durations of the previous jobs. Therefore, to facilitate the gradient computation in equation (\ref{derv}), we use Tensorflow 1.8 Python API to benefit from automatic gradient computation (automatic differentiation) \cite{abadi2016tensorflow}.

\section{Experiments} \label{section5}
In this section, we present our results for a randomly generated dataset, and CT scan images from the Lung Image Database Consortium image collection (LIDC-IDRI) \cite{armato2015data} with randomly generated appointment durations. 

\subsection{Randomly Generated Data}
Inspired by medical appointment dataset \cite{noshow2016data} on Kaggle, we generated a random dataset with $4$ categorical covariates (gender, day of month, time of day, intensity) and the response variable shows the duration of the appointment (non-negative). We assume the ground truth distribution of the duration times can be from multiple different distributions, but in all of them mean and standard deviation is fixed:
\begin{align*}
& \mu_{(x_2,x_3,x_4)} = 1+0.1x_2+0.4x_3+1.5x_4\\
& \sigma_{(x_1,x_4)} = 0.1+0.2x_1+0.2x_4
\end{align*}

Let us consider $4$ possible distributions for the response in here:
\begin{center}
\begin{tabular}{ll} \label{distsDef}
\textbf{Distribution} &\textbf{Description}\\
\hline
Normal (trunc.) & $\mathcal{N}\left(\mu_{(.)},\sigma_{(.)}\right)$\\
Logistic (trunc.) &  $LOG\left(\mu_{(.)},\sigma_{(.)}\right)$\\
Beta (scaled) &  $B\left(\alpha(\mu_{(.)},\sigma_{(.)}),\beta(\mu_{(.)},\sigma_{(.)})\right)$\\
Uniform  &  $\mathcal{U}\left(a(\mu_{(.)},\sigma_{(.)}),b(\mu_{(.)},\sigma_{(.)})\right)$
\end{tabular}
\end{center}
Later, we add some noise to the labels from $\mathcal{U}(-1,1)$. The appointment durations are not i.i.d, however, i.i.d is a conventional assumption for the inputs of the classic models in stochastic optimization.

To make sure predictions are non-negative, we scale the appoinment durations between $0$ and $1$. This transformation is necessary since our predictions are outputs of a sigmoid function. Since the covariates are categorical, we consider One-Hot encoding transformations of them for the training.

\paragraph{Architecture:} We considered a Fully-Connected (FC) architecture with two dense layers for this experiment. The input layer has $2+30+5+4=41$ nodes. The number of hidden nodes are chosen among powers of 2 (i.e. $2^n, 2 \leq n \leq 8$). Using grid search on the hyperparameters, we used the best performance among all of our runs. The weights are initialized using the Xavier initializer \cite{glorot2010understanding}. The output layer has only one node and returns the scaled predictions. For the choice of optimizer, we compared the results between SGD, ADAM \cite{kingma2014adam}, L-BFGS \cite{nocedal1980updating} and chose the best one in each run. The learning rate is also chosen from $[0.001, 0.1]$ range using a grid search. We used different strategies such as constant, time inverse decaying, and exponential decaying to change the learning rates in each epoch and for each run.

\begin{figure}
\centerline{\includegraphics[width=0.3\textwidth]{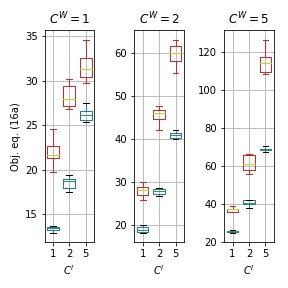}}
\caption{The Box-plots of the objective function (\ref{saaobj}) for different combinations of ($C^W$,$C^I$)}
\label{fig:saaObj}
\centering
\end{figure}

\subsection{Solutions of the Classic ASP}
The Maximum Likelihood Estimates (MLEs) of the mean and standard deviation parameters from the generated data are computed and used in this section. From the dataset, we have estimated $\hat{\mu}_{MLE}=5.9254$, and $\hat{\sigma}_{MLE}=2.2151$ for the truncated Normal distribution.

We use the Gurobi $7.0$ solver API for Python to solve the models. The solver parameters for tuning are all set to their default values. We use SAA with $10$ reps and equal number of samples to the original dataset for a model with $n=5$ (scheduling $5$ jobs on a single server). The Box-plots of the objective function (\ref{saaobj}) for different combinations of $C^W$ and $C^I$ parameters from the list of $[1,2,5]$ are depicted in Figure \ref{fig:saaObj}. The blue plots are related to the SAA model and the red ones are related to the DRO model.

\begin{figure}
\centerline{\includegraphics[width=0.4\textwidth]{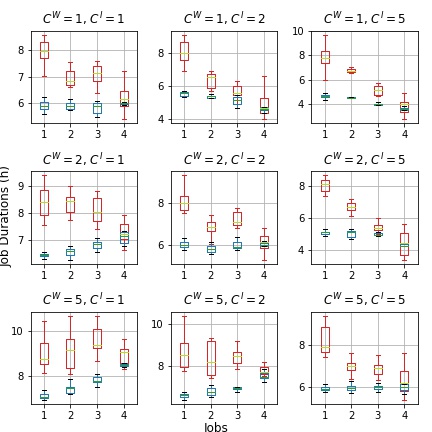}}
\caption{The Box-plots of the ASP solutions for different combinations of ($C^W$,$C^I$)} 
\label{fig:saaSol}
\centering
\end{figure}

As we can see in Figure \ref{fig:saaObj}, the base cost of the planning for $C^W=1$ and $C^I=1$ is estimated to be $13.40$ for SAA and $22.02$ for DRO and it is increasing proportionally by $C^W$ and $C^I$. Since DRO considers the worst-case distribution given the fixed mean and standard deviation, in all cases the total objective of DRO is higher than SAA. The actual decisions and their Box-plots are shown in Figure \ref{fig:saaSol}.

The interesting trends in Figure \ref{fig:saaSol} are related to the optimal scheduled appointments when changing $C^W$ and $C^I$. As we can see, in settings with high ratio of $\frac{C^W}{C^I}$, the length of the allowances are increasing by the order of the jobs and in settings with with high ratio of $\frac{C^I}{C^W}$, the allowances are decreasing by the order of the jobs. Job duration estimates from DRO are generally higher than SAA. This is due to the fact that DRO considers the worst-case distribution of the job durations whereas SAA assumes the exact distribution. Further, DRO is more sensitive to the ratio of the wait/idle costs.

\subsection{Comparison of IEO and SEO}
To make a fair comparison, we use Maximum Absolute Deviation (MAD) for the loss function of the SEO and the proposed loss function (\ref{ASPloss}) for the IEO. Moreover, we evaluate the actual objective value (\ref{saaobj}) for the solutions from both method. We also make sure all other parameters and hyperparameters of the networks are the same in both settings.

\begin{table}
\setlength{\tabcolsep}{2.5pt}
  \centering
  \caption{The objective function (\ref{saaobj}) for different combinations of ($C^W$,$C^I$)}
    \begin{tabular}{cc|ccc|ccc}
    \toprule
          & \textbf{} & \multicolumn{3}{c}{\textbf{SEO}} & \multicolumn{3}{c}{\textbf{IEO}} \\
    \multicolumn{1}{l}{\textbf{Dist.}} &  $\bm{C^w/C^I}$ & $1$  & $2$  & $5$  & $1$  & $2$  & $5$ \\
    \midrule
    \multirow{3}[6]{*}{\textbf{Norm.}} & $1$  & 10.56 & 15.63 & 30.86 & 10.86 & 14.92 & 19.93 \\
       & $2$  &  16.04 & 21.11 & 36.34 & 14.85 & 21.73 & 32.49 \\
          & $5$  & 32.48 & 37.55 & 52.78 & 20.64 & 32.40 & 54.44 \\
    \midrule
    \multirow{3}[6]{*}{\textbf{Log.}} & $1$  & 11.09 & 16.81 & 33.98 & 11.25 & 16.56 & 23.10 \\
          & $2$   & 16.46 & 22.18 & 39.35 & 15.81 & 23.79 & 36.54 \\
          & $5$  & 32.55 & 38.28 & 55.45 & 22.34 & 36.69 & 59.73  \\
    \midrule
    \multirow{3}[6]{*}{\textbf{Beta}} & $1$  & 10.61 & 15.99 & 31.53 & 10.94 & 14.87 & 20.23 \\
          & $2$  & 15.84 & 21.22 & 37.35 & 14.90 & 21.78 & 32.15 \\
          & $5$ & 31.53 & 36.91 & 53.04 & 20.40 & 32.48 & 54.70 \\
    \midrule
    \multirow{3}[6]{*}{\textbf{Uni.}} & $1$  &  10.50 & 15.81 & 31.72 & 10.91 & 14.75 & 20.16 \\
          & $2$   & 15.70 & 21.01 & 36.92 & 14.99 & 21.79 & 32.28 \\
          & $5$  & 31.30 & 36.61 & 52.52 & 20.42 & 32.28 & 54.49 \\
    \end{tabular}%
  \label{tab:tableSEOIEO}%
\end{table}%

As we can see in Table \ref{tab:tableSEOIEO}, the performance of the IEO is better than SEO in most of the settings. The most interesting part of this table is the objective function for the base parameters $C^W=1$ and $C^I=1$. In here, the solutions from the IEO is $10.56$, and from SEO to $10.86$ and both of them are better than the objective from classic model (13.40). This fact indicates that the i.i.d assumption is not valid.

From Table\ref{tab:tableSEOIEO}, we can also conclude that the choice of the distribution does not really affect the performance of the both models. As we expected, the performance of the SEO drops for the cases with high values of $C^W$ or $C^I$. This is due to the fact that in SEO we ignore these costs when we are estimating the job durations. The results clearly show how SEO can lead to sub-optimal solution.

\subsection{LIDC-IDRI Dataset}
The benefits of using DNNs for the contextual ASP do not end with bringing in context in decision making and creating a task-based loss function. In some cases, other ML algorithms such as SVM can be used, but they are often constrained by the assumption that we can define all the covariates. However, this is not necessarily true in many cases and making assumptions about non-numeric features can lead to low sensitivity. Deep learning is a fast-growing field that could be an ideal solution for these cases. This is due to the fact that we can learn features from raw images. In screening process for lung cancer, millions of CT scans will have to be analyzed. This is a controversial task even among radiologists. Further, radiologists often have to look through large volumes of these images that can lead to mistakes. Here, we skip the details about the screening process, and the dataset. We refer the interested readers to \cite{armato2011lung} for more details about lung nodule analysis. 

Consider a scenario in which we try to schedule patients for a second screening based on their initial phase results. Therefore, if the radiologist in the initial phase detects a nodule, the patient needs more time for the second phase. The process of scheduling patients in this manner can be an overwhelming burden for the radiologists. In this paper, we use the LIDC-IDRI dataset with 8723 CT scans. Similar to the previous example, we add more covariates (gender, day, hour) and we sample the response variable (appointment duration) from a Normal distribution using similar means and standard deviations to the previous experiment. The only difference in this experiment with the previous one is that the intensity is not a categorical column anymore and it is replaced by the image.

\paragraph{Architecture:} We used $3$ convolutional layers with $32$, $64$, and $128$ filters of $2\times2$ kernels followed by $2$ dense layers. The last dense layer has one output node with  a sigmoid activation corresponding to positive and negative cases. Using small kernel sizes and not using  max-pooling layers allow us to detect small areas of interest in the input images. We first resized all images to $55\times55$ pixels and then used several data augmentation methods including random cropping of $50\times50$, random horizontal flipping, color jittering. In addition to our surrogate loss function, we used a binary crossentropy loss together with Adam optimizer for the classification task.

The training process with the convolutional layers is significantly more expensive, therefore, we allow for more training epochs in this experiment. The accuracy of the  convolutional layers for the classification task increased up to $\%96$ after the augmentation. The results in Table \ref{tab:tableCNN} suggests that the IEO outperforms SEO in all cases, specifically for higher values of $C^W$ and $C^I$.
\begin{table}
\setlength{\tabcolsep}{2.5pt}
  \centering
  \caption{The objective function (\ref{saaobj}) for different combinations of ($C^W$,$C^I$)}
    \begin{tabular}{cc|ccc|ccc}
    \toprule
          & \textbf{} & \multicolumn{3}{c}{\textbf{IEO}} & \multicolumn{3}{c}{\textbf{SEO}} \\
    \multicolumn{1}{l}{\textbf{Dist.}} &  $\bm{C^w/C^I}$ & $1$  & $2$  & $5$  & $1$  & $2$  & $5$ \\
    \midrule
    \multirow{3}[6]{*}{\textbf{Norm.}} & $1$ & 11.88 & 16.01 & 21.34 & 11.90 & 19.15 & 40.92 \\
    & $2$ & 16.84 & 23.65 & 34.19 & 16.54 & 23.79 & 45.56 \\
    & $5$ & 23.89 & 24.10 & 59.37 & 30.47 & 37.72 & 59.48 \\
   \midrule
    \end{tabular}%
  \label{tab:tableCNN}%
\end{table}%

\section{Conclusions} \label{section6}
In this paper, we consider the contextual Appointment Scheduling Problem (ASP). If the probability distribution of the job durations is known for every possible combination of the data features, there is an exact solution for this problem. However, approximating a probability distribution is inaccurate, and impossible when data is scarce. The SEO approach to solve ASP is also shown to yield sub-optimal decisions. In particular, in cases that the idling or waiting cost is significantly higher than the other cost. To address this issue, we proposed the IEO approach with a surrogate loss function in which we estimate the durations using non-linear models. This approach does not require the knowledge of the job durations probability distribution and uses only historical data. We have shown the validity of our approach though the comparisons of the solutions from IEO with other approaches.
\bibliography{main.bib}

\begin{thebibliography}{}

\bibitem[\protect\citeauthoryear{Abadi \bgroup et al\mbox.\egroup
  }{2016}]{abadi2016tensorflow}
Abadi, M.; Barham, P.; Chen, J.; Chen, Z.; Davis, A.; Dean, J.; Devin, M.;
  Ghemawat, S.; Irving, G.; Isard, M.; et~al.
\newblock 2016.
\newblock Tensorflow: A system for large-scale machine learning.
\newblock In {\em OSDI}, volume~16,  265--283.

\bibitem[\protect\citeauthoryear{Armato~III \bgroup et al\mbox.\egroup
  }{2011}]{armato2011lung}
Armato~III, S.~G.; McLennan, G.; Bidaut, L.; McNitt-Gray, M.~F.; Meyer, C.~R.;
  Reeves, A.~P.; Zhao, B.; Aberle, D.~R.; Henschke, C.~I.; Hoffman, E.~A.;
  et~al.
\newblock 2011.
\newblock The lung image database consortium (lidc) and image database resource
  initiative (idri): a completed reference database of lung nodules on ct
  scans.
\newblock {\em Medical physics} 38(2):915--931.

\bibitem[\protect\citeauthoryear{Armato~III \bgroup et al\mbox.\egroup
  }{2015}]{armato2015data}
Armato~III, S.~G.; McLennan, G.; Bidaut, L.; McNitt-Gray, M.~F.; Meyer, C.~R.;
  Reeves, A.~P.; and Clarke, L.~P.
\newblock 2015.
\newblock Data from lidc-idri. the cancer imaging archive.

\bibitem[\protect\citeauthoryear{Ban and Rudin}{2018}]{ban2017big}
Ban, G.-Y., and Rudin, C.
\newblock 2018.
\newblock The big data newsvendor: Practical insights from machine learning.
\newblock {\em Forthcoming in Operations Research}.

\bibitem[\protect\citeauthoryear{Bertsimas and
  Kallus}{2014}]{bertsimas2014predictive}
Bertsimas, D., and Kallus, N.
\newblock 2014.
\newblock From predictive to prescriptive analytics.
\newblock {\em arXiv preprint arXiv:1402.5481}.

\bibitem[\protect\citeauthoryear{Cardoen, Demeulemeester, and
  Beli{\"e}n}{2010}]{cardoen2010operating}
Cardoen, B.; Demeulemeester, E.; and Beli{\"e}n, J.
\newblock 2010.
\newblock Operating room planning and scheduling: A literature review.
\newblock {\em European journal of operational research} 201(3):921--932.

\bibitem[\protect\citeauthoryear{Cayirli, Veral, and
  Rosen}{2006}]{cayirli2006designing}
Cayirli, T.; Veral, E.; and Rosen, H.
\newblock 2006.
\newblock Designing appointment scheduling systems for ambulatory care
  services.
\newblock {\em Health care management science} 9(1):47--58.

\bibitem[\protect\citeauthoryear{Denton and Gupta}{2003}]{denton2003sequential}
Denton, B., and Gupta, D.
\newblock 2003.
\newblock A sequential bounding approach for optimal appointment scheduling.
\newblock {\em Iie Transactions} 35(11):1003--1016.

\bibitem[\protect\citeauthoryear{Donti, Amos, and Kolter}{2017}]{donti2017task}
Donti, P.~L.; Amos, B.; and Kolter, J.~Z.
\newblock 2017.
\newblock Task-based end-to-end model learning in stochastic optimization.
\newblock {\em arXiv preprint arXiv:1703.04529}.

\bibitem[\protect\citeauthoryear{Elmachtoub and
  Grigas}{2017}]{elmachtoub2017smart}
Elmachtoub, A.~N., and Grigas, P.
\newblock 2017.
\newblock Smart" predict, then optimize".
\newblock {\em arXiv preprint arXiv:1710.08005}.

\bibitem[\protect\citeauthoryear{Gallego and
  Moon}{1993}]{gallego1993distribution}
Gallego, G., and Moon, I.
\newblock 1993.
\newblock The distribution free newsboy problem: review and extensions.
\newblock {\em Journal of the Operational Research Society} 44(8):825--834.

\bibitem[\protect\citeauthoryear{Glorot and
  Bengio}{2010}]{glorot2010understanding}
Glorot, X., and Bengio, Y.
\newblock 2010.
\newblock Understanding the difficulty of training deep feedforward neural
  networks.
\newblock In {\em Proceedings of the thirteenth international conference on
  artificial intelligence and statistics},  249--256.

\bibitem[\protect\citeauthoryear{Gould \bgroup et al\mbox.\egroup
  }{2016}]{gould2016differentiating}
Gould, S.; Fernando, B.; Cherian, A.; Anderson, P.; Cruz, R.~S.; and Guo, E.
\newblock 2016.
\newblock On differentiating parameterized argmin and argmax problems with
  application to bi-level optimization.
\newblock {\em arXiv preprint arXiv:1607.05447}.

\bibitem[\protect\citeauthoryear{Gupta and Denton}{2008}]{gupta2008appointment}
Gupta, D., and Denton, B.
\newblock 2008.
\newblock Appointment scheduling in health care: Challenges and opportunities.
\newblock {\em IIE transactions} 40(9):800--819.

\bibitem[\protect\citeauthoryear{Gurvich, Luedtke, and
  Tezcan}{2010}]{gurvich2010staffing}
Gurvich, I.; Luedtke, J.; and Tezcan, T.
\newblock 2010.
\newblock Staffing call centers with uncertain demand forecasts: A
  chance-constrained optimization approach.
\newblock {\em Management Science} 56(7):1093--1115.

\bibitem[\protect\citeauthoryear{Hoppen}{2016}]{noshow2016data}
Hoppen, J.
\newblock 2016.
\newblock Data from kaggle datasets. the medical appointment no shows.

\bibitem[\protect\citeauthoryear{Kingma and Ba}{2014}]{kingma2014adam}
Kingma, D.~P., and Ba, J.
\newblock 2014.
\newblock Adam: A method for stochastic optimization.
\newblock {\em arXiv preprint arXiv:1412.6980}.

\bibitem[\protect\citeauthoryear{Kleywegt, Shapiro, and Homem-de
  Mello}{2002}]{kleywegt2002sample}
Kleywegt, A.~J.; Shapiro, A.; and Homem-de Mello, T.
\newblock 2002.
\newblock The sample average approximation method for stochastic discrete
  optimization.
\newblock {\em SIAM Journal on Optimization} 12(2):479--502.

\bibitem[\protect\citeauthoryear{Ko{\c{c}}a{\u{g}}a, Armony, and
  Ward}{2015}]{koccauga2015staffing}
Ko{\c{c}}a{\u{g}}a, Y.~L.; Armony, M.; and Ward, A.~R.
\newblock 2015.
\newblock Staffing call centers with uncertain arrival rates and co-sourcing.
\newblock {\em Production and Operations Management} 24(7):1101--1117.

\bibitem[\protect\citeauthoryear{Koenker and
  Bassett~Jr}{1978}]{koenker1978regression}
Koenker, R., and Bassett~Jr, G.
\newblock 1978.
\newblock Regression quantiles.
\newblock {\em Econometrica: journal of the Econometric Society}  33--50.

\bibitem[\protect\citeauthoryear{Lindley}{1952}]{lindley1952theory}
Lindley, D.~V.
\newblock 1952.
\newblock The theory of queues with a single server.
\newblock In {\em Mathematical Proceedings of the Cambridge Philosophical
  Society}, volume~48,  277--289.
\newblock Cambridge University Press.

\bibitem[\protect\citeauthoryear{Mak, Rong, and
  Zhang}{2014}]{mak2014appointment}
Mak, H.-Y.; Rong, Y.; and Zhang, J.
\newblock 2014.
\newblock Appointment scheduling with limited distributional information.
\newblock {\em Management Science} 61(2):316--334.

\bibitem[\protect\citeauthoryear{Mittal, Schulz, and
  Stiller}{2014}]{mittal2014robust}
Mittal, S.; Schulz, A.~S.; and Stiller, S.
\newblock 2014.
\newblock Robust appointment scheduling.
\newblock In {\em LIPIcs-Leibniz International Proceedings in Informatics},
  volume~28.
\newblock Schloss Dagstuhl-Leibniz-Zentrum fuer Informatik.

\bibitem[\protect\citeauthoryear{Nocedal}{1980}]{nocedal1980updating}
Nocedal, J.
\newblock 1980.
\newblock Updating quasi-newton matrices with limited storage.
\newblock {\em Mathematics of computation} 35(151):773--782.

\bibitem[\protect\citeauthoryear{Oroojlooyjadid, Snyder, and
  Takac}{2016}]{OroojlooyjadidS16}
Oroojlooyjadid, A.; Snyder, L.~V.; and Takac, M.
\newblock 2016.
\newblock Applying deep learning to the newsvendor problem.
\newblock {\em CoRR} abs/1607.02177.

\bibitem[\protect\citeauthoryear{Shapiro, Dentcheva, and
  Ruszczy{\'n}ski}{2009}]{shapiro2009lectures}
Shapiro, A.; Dentcheva, D.; and Ruszczy{\'n}ski, A.
\newblock 2009.
\newblock {\em Lectures on stochastic programming: modeling and theory}.
\newblock SIAM.

\bibitem[\protect\citeauthoryear{Song \bgroup et al\mbox.\egroup
  }{2016}]{song2016training}
Song, Y.; Schwing, A.; Urtasun, R.; et~al.
\newblock 2016.
\newblock Training deep neural networks via direct loss minimization.
\newblock In {\em International Conference on Machine Learning},  2169--2177.

\bibitem[\protect\citeauthoryear{Taylor}{2000}]{taylor2000quantile}
Taylor, J.~W.
\newblock 2000.
\newblock A quantile regression neural network approach to estimating the
  conditional density of multiperiod returns.
\newblock {\em Journal of Forecasting} 19(4):299--311.

\bibitem[\protect\citeauthoryear{Weiss}{1990}]{weiss1990models}
Weiss, E.~N.
\newblock 1990.
\newblock Models for determining estimated start times and case orderings in
  hospital operating rooms.
\newblock {\em IIE transactions} 22(2):143--150.

\bibitem[\protect\citeauthoryear{White}{1992}]{white1992nonparametric}
White, H.
\newblock 1992.
\newblock Nonparametric estimation of conditional quantiles using neural
  networks.
\newblock In {\em Computing Science and Statistics}. Springer.
\newblock  190--199.

\bibitem[\protect\citeauthoryear{Zhang and Gao}{2017}]{zhang2017assessing}
Zhang, Y., and Gao, J.
\newblock 2017.
\newblock Assessing the performance of deep learning algorithms for newsvendor
  problem.
\newblock In {\em International Conference on Neural Information Processing},
  912--921.
\newblock Springer.

\bibitem[\protect\citeauthoryear{Zhang, Shen, and
  Erdogan}{2017}]{zhang2017distributionally}
Zhang, Y.; Shen, S.; and Erdogan, S.~A.
\newblock 2017.
\newblock Distributionally robust appointment scheduling with moment-based
  ambiguity set.
\newblock {\em Operations Research Letters} 45(2):139--144.

\end{thebibliography}
\bibliographystyle{aaai}
\end{document}